\documentclass[11pt]{amsart}
\usepackage{amsmath,amsthm, amscd, amssymb, amsfonts}

\newcommand{\we}{{\widetilde e}}

\newcommand{\ka}{{\mathfrak k}}
\newcommand{\ko}{{\mathfrak k}_0}

\newcommand{\p}{{\mathfrak p}}
\newcommand{\po}{{\mathfrak p}_0}
\newcommand{\tc}{{\mathfrak t}_0}
\newcommand{\ao}{{\mathfrak a}_0}
\newcommand{\ag}{\mathfrak a}
\newcommand\dc{\operatorname{dc}}
\newcommand\dnc{\operatorname{dnc}}

\newcommand\diag{\operatorname{diag}}

\newcommand\re{\operatorname{Re}}

\newcommand{\dynkin}{\mu}
\newcommand{\tdynkin}{\varsigma_{\mu}}
\newcommand{\tid}{\varsigma}

\newcommand{\jdynkin}{\omega_{\mu, J}}
\newcommand{\jid}{\omega_{J}}

\newcommand{\comp}{\omega}

\newcommand{\go}{{\mathfrak g}_0}
\newcommand{\g}{{\mathfrak g}}
\newcommand{\ho}{{\mathfrak h}_0}
\newcommand{\h}{{\mathfrak h}}

\newcommand{\ro}{{r}_0}
\newcommand{\rom}{{r}_\Omega}
\newcommand{\rla}{{r}_\Lambda}

\newcommand{\lgot}{{\mathfrak l}}
\newcommand{\lgo}{{\mathfrak l}_0}
\newcommand{\lgor}{\lgot^{\R}}
\newcommand{\invo}{{\sigma}}

\newcommand\Aut{\operatorname{Aut}}

\newcommand{\Z}{{\mathbb Z}}

\newcommand{\R}{{\mathbb R}}
\newcommand{\Cc}{{\mathbb C}}

\numberwithin{equation}{section}

\newcommand{\id}{\mathop{\rm id\,}}

\newcommand{\ad}{\mbox{\rm ad\,}}
\newcommand{\cent}{\mathop{\rm Cent\,}}

\newcommand{\card}{\mbox{\rm card\,}}

\theoremstyle{plain}

\newtheorem{theorem}{Theorem}
\newtheorem{lema}{Lemma}[section]

\newtheorem{prop}[lema]{Proposition}

\theoremstyle{definition}
\newtheorem{definition}[lema]{Definition}

\theoremstyle{remark}
\newtheorem{obs}[lema]{Remark}

\def\pf{\begin{proof}}
\def\epf{\end{proof}}

\theoremstyle{remark}

\begin{document}

\renewcommand{\baselinestretch}{1.2}
\renewcommand{\thefootnote}{}
\thispagestyle{empty}
%\vspace*{2in}
\title[Simple real Lie bialgebras]{On simple real Lie bialgebras}
\author[Andruskiewitsch and Jancsa]{ Nicol\'as Andruskiewitsch and
Patricia Jancsa}
\address{Facultad de Matem\'atica, Astronom\'\i a y F\'\i sica,
Universidad Nacional de C\'ordoba
\newline
\indent CIEM -- CONICET
\newline
\indent (5000) Ciudad Universitaria, C\'ordoba, Argentina}
\email{andrus@mate.uncor.edu, \ jancsa@mate.uncor.edu}
\thanks{This work was partially supported by CONICET,
Agencia C\'ordoba Ciencia, ANPCyT  and Secyt (UNC). Results presented in this paper 
are part of the PhD-thesis of the second-named author, written under the advise 
of the first-named author at the University of C\'ordoba, Argentina.}
\date{Revised version of  March 6th, 2003. (First version as of October 14th, 2002, appeared as IHES Preprint
M/02/81).}
\begin{abstract}
The explicit list of all almost factorizable Lie bialgebra structures 
on real absolutely simple Lie algebras is given.
By ``absolutely simple" we mean a  real Lie algebra
whose complexification is  simple.
\end{abstract}
\maketitle

\section*{Introduction}

The theory of Poisson-Lie groups occupies a central place in the
theory of Poisson manifolds. The category of connected,
simply-connected Poisson-Lie groups is equivalent to the category
of Lie bialgebras \cite{D}. Therefore, a basic problem in the
theory of Poisson manifolds is the classification of Lie
bialgebras. A fundamental contribution to this question is the
Theorem of Belavin and Drinfeld \cite{BD}, which contains the
classification of all the simple factorizable complex Lie
bialgebras. 

A finite-dimensional real Lie algebra is \emph{absolutely simple} 
if its complexification is a simple complex Lie algebra; this 
notation agrees with the tradition in Lie theory to call 
``absolutely simple" an object that remains simple after any 
extension of scalars \cite{Di, T}. It is well-known that a 
simple real Lie algebra is either absolutely simple, or 
it is the realification of a complex simple Lie algebra.

In this paper, we obtain the following result.

\begin{theorem}\label{main1} Let $(\go, \delta)$ be an absolutely 
simple real Lie bialgebra.
Then either $(\go, \delta)$ is triangular, or else it is
determined by the pair $(\go, \ro)$ as in Table \ref{tablauno}, up
to isomorphisms of real Lie bialgebras.

\bigbreak Precisely, there exists a unique Cartan subalgebra $\h$
of the complexification $\g$ of $\go$; a system of simple roots
$\Delta \subset \Phi(\g, \h)$; an involution $\invo = \tdynkin$ or
$\jdynkin$  as in \eqref{tdynkin}, \eqref{jdynkin}; a BD-triple $(\Gamma_1,
\Gamma_2, T)$; $t\in \R_{>0} \cup
i\R_{>0}$; and a continuous parameter $\lambda \in \h^{\otimes
2}$; all these data subject to the restrictions in Table
\ref{tablauno}; such that $(\go, \delta)$ is isomorphic as real
Lie bialgebra to $(\g^{\invo},
\partial \ro)$, where

$$\ro =t\left( \frac 12 \sum_{\alpha, \beta \in \Delta}
 \lambda_{\alpha, \beta} h_{\alpha} \wedge h_\beta + \frac 12
\sum_{\alpha \in \Phi^+} x_{-\alpha} \wedge x_\alpha
        + \sum_{\alpha, \beta \in \Phi^+,  \alpha \prec 
\beta}
                               x_{-\alpha} \wedge x_\beta\right).$$
Here, $\lambda - \lambda^{21} = \sum_{\alpha, \beta \in
\Delta} \lambda_{\alpha, \beta} h_{\alpha} \wedge
h_\beta$, and $h_{\alpha}$, $x_\beta$ have the usual meaning, see
subsections \ref{conv}, \ref{inv}.

Two data in the table give rise to isomorphic Lie bialgebras if and only 
if they belong to the same line and the rest of the data is conjugated by 
an automorphism of the Dynkin diagram of order 2.
\end{theorem}

\begin{table}[t]
\begin{center}
\begin{tabular}{|p{3,5cm}|p{2,2cm}|p{1cm}|p{2,1cm}|p{3,5cm}|p{0,6cm}|}
\hline 

{\bf $\go$} & $\invo$ & {\bf Type} &
{\bf BD-triple} & {\bf Continuous \newline parameter} & $t\in $

\\ \hline $\g_{\R}$& $\tid$ & All.  & All. & $\lambda_{\alpha, \beta} \in \R$  &
$\R$

\\ \hline $\mathfrak{su}(n, n  + 1)$ \newline
$\mathfrak{su}(n + 1, n  + 1)$ \newline
$\mathfrak{so}(n - 1, n  + 1)$\newline
EII
 & $\tdynkin$, $\dynkin \neq  \id$
&  $A_{2n}$ \newline $A_{2n + 1}$\newline $D_{n}$
\newline $E_{6}$ & $\dynkin$-stable
& ${}^{\,}$\newline
$\lambda_{\alpha, \beta} = \overline{\lambda_{\dynkin(\alpha), 
\dynkin(\beta)}}$ & $\R$

\\ \hline $\g_c$ & $\comp$  & All.  & $\Gamma_1 = \Gamma_2 = 
\emptyset$
& $\lambda_{\alpha, \beta} \in i\R$   & $i\R$

\\ \hline  $\mathfrak{su}(j, n + 1 - j)$
\newline $\mathfrak{so}(2j, 2(n -j) + 1)$
\newline $\mathfrak{sp}(j, n  - j)$
\newline $\mathfrak{so}(2j, 2n - 2j)$
\newline $\mathfrak{so}^*(2n)$
\newline EII,  EIII
\newline EV, EVI, EVII
\newline EVIII, EIX
\newline FI, FII
\newline G
& $\jid$ &  $A_n$
\newline $B_n$
\newline $C_n$
\newline $D_n$
\newline $D_n$
\newline $E_6$
\newline $E_7$
\newline $E_8$
\newline $F_4$
\newline $G_2$
& $\Gamma_1 = \Gamma_2 = \emptyset$ & $\lambda_{\alpha, \beta}
\in i\R$   & $i\R$

\\ \hline  $\mathfrak{sl}(n+1, \R)$
\newline $\mathfrak{sl}(\frac{n+1}2, \mathbb H)$ \newline  
$\mathfrak{so}(1, 2n  - 1)$
\newline  $\mathfrak{so}(2j +1, 2(n - j)  - 1)$ \newline EI, EIV
& $\jdynkin$, $\dynkin \neq  \id$ & 
$A_{n}$
\newline $A_{2n + 1}$
\newline $D_{n}$ \newline $D_{n}$ \newline $E_{6}$
& $\dynkin$-antistable
& ${}^{\,}$\newline $\lambda_{\alpha, \beta} = 
-\overline{\lambda_{\dynkin(\alpha), \dynkin(\beta)}}$ &
$i\R$
\\ \hline
\end{tabular}
\end{center}

\

\emph{Explanation of the table.} We recall the notions of BD-triples
and continuous parameter in subsection \ref{BDsection}; see
Definition \ref{dynkinstable} for the notion of $\dynkin$-stable
or antistable. See also Remark \ref{isom} for a discussion on the involution. 
The notation for exceptional $\go$'s is standard
and goes back to \'E. Cartan \cite{cartan}.

\

\caption{Classification of almost-factorizable simple real Lie 
bialgebras}\label{tablauno}

\end{table}

\bigbreak Some considerations about real simple Lie bialgebras are
already present in \cite{LQ, A, CGR, Ch, KRR}. The compact case is
well-known, see for example \cite{KS, LW, M}. All real classical
$r$-matrices arising from semisimple Lie algebras are listed in
\cite{CGR}. A. Panov described in \cite{P1}, \cite{P2}, all
posible Manin triples in a suitable class up to gauge and weak
equivalence; see also \cite{De}. One might be able to infer our
main result from these papers, with extra work. However, we do not
follow this path and deduce the classification directly from the
Theorem of Belavin and Drinfeld. The key step in our argument is
Lemma \ref{necesario}, where we attach a Cartan subalgebra and a
choice of positive roots to the pair $(\go, \ro)$.

\bigbreak The complexifications of the Lie bialgebras associated
to the pairs $(\go, \ro)$ as in Table \ref{tablauno} are
factorizable, as defined in \cite{RS}. This is not the case for
the real Lie bialgebras themselves; some of them are
quasitriangular, some of them are not. We use the term ``almost
factorizable" to refer to Lie bialgebras that are factorizable after 
extension of scalars.

\bigbreak
The organization of this paper is as follows. In the first section we 
recall
several well-known facts about real Lie algebras (subsections 1.1 and 
1.2)
and  Lie bialgebras (subsection 1.3). The material on real Lie algebras 
is standard
but we decided to present it for completeness of
the arguments leading to the main result. In the second section
we prove two lemmas that imply the main result.
In the last section we compute the Drinfeld double and the dual Lie bialgebra 
of a real absolutely simple Lie bialgebra. Our result is a consequence
of the determination of the Drinfeld double of a factorizable Lie bialgebra \cite{RS}.

\section*{Acknowledgements}
We thank A. Tiraboschi for many conversations, and the referee for some interesting 
remarks and suggestions. The first author
expresses his gratitude to IHES, where the final part of the work
in this article was done.

\section{Preliminaries}

\subsection{Conventions}\label{conv}
All Lie algebras are finite-dimensional, unless explicitly stated.
We shall say that a Lie bialgebra is simple if the underlying Lie 
algebra is simple.

Let $\go$ be a semisimple real Lie algebra and let $\g$ be its
complexification. The Killing form on $\g$, or $\go$, is denoted
by $(\quad\vert\quad)$. The Casimir element of $\g$ is denoted by
$\Omega \in \g \otimes \g$; that is, $\Omega = \sum x_i \otimes
x^i$ where $(x_i)$, $(x^i)$ is any pair of dual basis with respect
to the Killing form of $\g$. Clearly $\Omega \in \go \otimes \go$.

Let $\h\subset \g$ be a Cartan subalgebra; we denote by $\Phi =
\Phi(\g, \h)$ the set of roots. If $\Delta$ is a system of simple
roots, we denote by $\Phi^+$ the corresponding set of positive
roots. We denote by $\Omega_{0}$ the component in $\h \otimes \h$
of $\Omega$; that is, $\Omega _{0}= \sum h_i \otimes h^i$ where
$(h_i)$, $(h^i)$ is any pair of dual basis in $\h$ with respect to
the restriction of the Killing form of $\g$ to $\h$. If $\alpha
\in \h^*$, we denote by $h_{\alpha} \in \h$ the unique element in
$\h$ such that $\alpha(H) = (h_{\alpha} \vert H)$ for all $H \in
\h$.

\bigbreak
If $V$ is a complex vector space and $\invo: V \to V$ is a sesquilinear 
map
with $\invo^2 = \id$, then $\invo^*: V^* \to V^*$ denotes the 
sesquilinear map
with $\invo^*(\alpha) (h) = \overline{\alpha(\invo(h))}$, $\alpha \in 
V^*$,
$h\in V$.

\subsection{Involutions}\label{inv}

Let $\g$ be a complex simple Lie algebra, $\h$ a Cartan subalgebra
of $\g$ , $\Delta \subset \Phi(\g, \h)$ a set of simple roots.
Let us choose elements $e_{\alpha}\in \g_{\alpha}$,
$e_{-\alpha}\in \g_{-\alpha}$, $\alpha \in \Delta$, such that
$[e_{\alpha}, e_{-\alpha}] = h_{\alpha}$.

\bigbreak Let $\dynkin: \Delta \to \Delta$ be an automorphism of
the Dynkin diagram of order 1 or 2. Let $J$ be any subset of the
set $\Delta^{\dynkin}$ of roots fixed by $\dynkin$; let $\chi_J:
\Delta \to \{0,1\}$ be the characteristic function of $J$ . Then
there exist unique sesquilinear Lie algebra involutions
$\tdynkin$, $\jdynkin$ of $\g$ given respectively by
\begin{alignat}{2}
\label{tdynkin} \tdynkin (e_{\alpha}) &= e_{\dynkin(\alpha)},
&\qquad \tdynkin (e_{-\alpha}) &= e_{-\dynkin(\alpha)}, \\
\label{jdynkin} \jdynkin (e_{\alpha}) &= (-1)^{\chi_J(\alpha)}
e_{-\dynkin(\alpha)}, &\qquad \jdynkin (e_{-\alpha}) &=
(-1)^{\chi_J(\alpha)} e_{\dynkin(\alpha)},
\end{alignat}
for all $\alpha \in \Delta$. This follows at once from Serre's
theorem on the presentation of $\g$, see for instance \cite{Kn}.
Necessarily, $\tdynkin (h_{\alpha}) = h_{\dynkin(\alpha)}$,
$\jdynkin (h_{\alpha}) = -h_{\dynkin(\alpha)}$. We shall
abbreviate $$\tid = \tid_{\id}, \qquad \jid = \comp_{\id, J},
\qquad \comp = \comp_{\Delta}.$$ Thus $\comp$ is the Chevalley
involution of $\g$, with respect to $\h$ and $\Delta$,  and the
fixed point set of $\comp$ is a compact form of $\g$, denoted by
$\g_c$.

\begin{lema}\label{involuciones}
Let $\invo: \g \to \g$ be a sesquilinear Lie algebra
involution such that $\invo (\h) = \h$. Then $\invo^* (\Phi) =
\Phi$ and $\invo(\g_{\alpha}) = \g_{\invo^* (\alpha)}$.

\medbreak (a). If $\invo^* (\Delta) = \Delta$, there exists a
choice of elements $e_{\alpha}\in \g_{\alpha}$, $e_{-\alpha}\in
\g_{-\alpha}$, $\alpha \in \Delta$ as above such that $\invo =
\tdynkin$ for a unique automorphism $\dynkin: \Delta \to \Delta$
of the Dynkin diagram of order 1 or 2.

\medbreak (b). If $\invo (\Delta)^* = -\Delta$, there exists a
choice of elements $e_{\alpha}\in \g_{\alpha}$, $e_{-\alpha}\in
\g_{-\alpha}$, $\alpha \in \Delta$ as above such that $\invo =
\jdynkin$ for a unique automorphism $\dynkin: \Delta \to \Delta$

of the Dynkin diagram of order 1 or 2, and a unique subset $J$ of
$\Delta^{\dynkin}$.
\end{lema}

\pf It is well-known that $(x\vert y) = \overline{(\invo(x) \vert
\invo(y))}$ for all $x$, $y\in \g$, see for example \cite[p.
180]{H}. If $H \in \h$ and $\alpha \in \Phi$, then for all $X\in
\g_{\alpha}$, one has $$[H, \invo(X)] = \invo([\invo(H), X]) =
\overline{\alpha(\invo(H))} \invo(X) = \invo^*(\alpha)(H)
\invo(X).$$ Hence $\invo^* (\Phi) = \Phi$ and
$\invo(\g_{\alpha}) = \g_{\invo^* (\alpha)}$.

Assume that $\invo^* (\Delta) = \pm \Delta$. Let $(a_{\alpha,
\beta})_{\alpha, \beta\in \Delta}$ be the Cartan matrix of $\g$.
Let $\dynkin: \Delta \to \Delta$ be given by $\dynkin = \pm
\invo^*$, according to the case. Then $$a_{\alpha, \beta} =
2\frac{(\alpha\vert \beta)}{(\beta \vert \beta)} =
\overline{{2\frac{(\invo^*(\alpha) \vert
\invo^*(\beta))}{(\invo^*(\beta) \vert \invo^*(\beta))}}^{\ }}
=a_{\dynkin(\alpha), \dynkin(\beta)},$$ hence $\dynkin$ is an
automorphism of the Dynkin diagram, and clearly it has order 1 or
2.

Let $e_{\alpha}\in \g_{\alpha}$, $e_{-\alpha}\in \g_{-\alpha}$,
$\alpha \in \Delta$ be as above. Let $c_{\alpha} \in \Cc -0$ be

such that $\invo(e_{\alpha}) = c_{\alpha} e_{\invo^*(\alpha)}$,
$\alpha \in \pm \Delta$. Then, for all $\alpha \in \Delta$, we
have
\begin{equation}\label{normuno} c_{\alpha} c_{-\alpha} = 1,
\end{equation} since $1 = (e_{\alpha}\vert e_{-\alpha}) =
\overline{(\invo(e_{\alpha}) \vert \invo(e_{-\alpha}))} =
\overline{c_{\alpha} c_{-\alpha}(e_{\invo^*(\alpha)}) \vert
e_{-\invo^*(\alpha)})}= \overline{c_{\alpha} c_{-\alpha}}$. Also,
$\invo^2 = \id$ implies that
\begin{equation}\label{normdos} c_{\alpha} \overline{c_{\invo^*\alpha}} 
= 1,
\end{equation} for all $\alpha \in \pm \Delta$.

\bigbreak We prove (a). We want to replace $e_{\alpha}$,
$e_{-\alpha}$ by $\we_{\alpha} = d_{\alpha}e_{\alpha}$,
$\we_{-\alpha} = (d_{\alpha})^{-1}e_{-\alpha}$, for well-chosen
non-zero scalars $d_{\alpha}$, $\alpha \in \Delta$, such that $$
\invo (\we_{\alpha}) = \we_{\dynkin(\alpha)}, \qquad \invo
(\we_{-\alpha}) = \we_{-\dynkin(\alpha)}. $$ But $\invo
(\we_{\alpha}) = \overline{d_{\alpha}}c_{\alpha}
(d_{\dynkin(\alpha)})^{-1} \we_{\dynkin(\alpha)}$ and $\invo
(\we_{-\alpha}) = \overline{(d_{\alpha})^{-1}}c_{-\alpha}
(d_{\dynkin(\alpha)}) \we_{-\dynkin(\alpha)}$, thus we need to
find $d_{\alpha}$ such that $ c_{\alpha} = \overline{ (
d_{\alpha})^{-1}} d_{\dynkin(\alpha)}$. The existence of such
$d_{\alpha}$ is clear, and the uniqueness of $\dynkin$ is evident;
(a) follows.

\bigbreak We prove (b). We want to replace $e_{\alpha}$,
$e_{-\alpha}$ by $\we_{\alpha} = d_{\alpha}e_{\alpha}$,
$\we_{-\alpha} = (d_{\alpha})^{-1}e_{-\alpha}$, for well-chosen
non-zero scalars $d_{\alpha}$, $\alpha \in \Delta$. We have
$\invo (\we_{\alpha}) = \overline{d_{\alpha}}c_{\alpha}
e_{-\dynkin(\alpha)} = \overline{d_{\alpha}}c_{\alpha}
d_{\dynkin(\alpha)} \we_{-\dynkin(\alpha)}$ and $\invo
(\we_{-\alpha}) = \overline{(d_{\alpha})^{-1}}c_{-\alpha}
(d_{\dynkin(\alpha)})^{-1} \we_{\dynkin(\alpha)}$.

\bigbreak Assume that $\dynkin(\alpha) \neq \alpha$; choose
$d_{\alpha}$, $d_{\dynkin(\alpha)}$ such that $c_{\alpha} =
(\overline{d_{\alpha}} d_{\dynkin(\alpha)})^{-1}$. Then $\invo
(\we_{\alpha}) =  \we_{-\dynkin(\alpha)}$, but also $\invo
(\we_{-\alpha}) =  \we_{\dynkin(\alpha)}$ by \eqref{normuno} and
$\invo (\we_{\pm \dynkin(\alpha)}) =  \we_{\mp \alpha}$
by \eqref{normuno} and \eqref{normdos}.

\bigbreak Assume that $\dynkin(\alpha) = \alpha$. Then
$c_{\alpha}\in \R$ by \eqref{normuno} and \eqref{normdos}. Let $J
:= \{\alpha\in \Delta: c_{\alpha} < 0\}$. Clearly, we can choose
$d_{\alpha}$ such that $c_{\alpha}  \vert d_{\alpha}\vert^2 = -1$,
resp. 1, if $\alpha\in J$, resp. if $\alpha\notin J$.

The uniqueness of $\dynkin$ and $J$ is evident, and (b) follows.
\epf

\begin{obs}\label{comp-inv}
The change of generators
$e_{\alpha}$, $e_{-\alpha}$ by $\we_{\alpha} =
d_{\alpha}e_{\alpha}$, $\we_{-\alpha} =
(d_{\alpha})^{-1}e_{-\alpha}$ amounts to composing with an
automorphism of $\g$ which preserves $\h$.
\end{obs}

\subsection{Real Lie algebras and Vogan diagrams}
We first briefly recall the main ingredients of the theory of real 
simple Lie
algebras, in particular Vogan diagrams; see \cite[Chapter VI]{Kn}. Next 
we apply these methods
to the study of the real Lie algebras corresponding to the involutions
in the preceding subsection.

\bigbreak
Let $\g$ be a complex simple Lie algebra and let $\invo$ be a
sesquilinear Lie algebra involution of $\g$.
Let $\go$ be the real form $\go = \g^{\invo}$ of $\g$.
Let $\theta_0$ be a Cartan
involution of $\go$ and let $\ho$ be a Cartan  subalgebra of $\go$ with
$\theta(\ho)= \ho$.
Let $\theta: \g \to \g$ be the complexification of $\theta_0$.
Let $\go = \ko \oplus \po$, respectively $\g = \ka \oplus \p$,  be the 
Cartan decomposition
associated to $\theta_0$, respectively $\theta$.
Denote $\tc := \ho \cap \ko$, $\ao := \ho \cap \po$.
Let $\dc(\ho) = \dim (\tc)$,
$\dnc(\ho) = \dim (\ao) = \dim \ho - \dc(\ho)$.
The Cartan subalgebra $\ho$ is \emph{maximally compact}, resp.
\emph{maximally non-compact} if $\dc(\ho)$ is maximal, resp. minimal, 
among
the $\dc$ of $\theta_0$-stable Cartan subalgebras of $\go$.
The \emph{rank} of the associated symmetric space coincides with
$\dnc(\ho)$ if $\ho$ is a maximally non-compact Cartan subalgebra.

Recall that a root $\alpha \in \Phi$ is called

\begin{itemize}
\item \emph{imaginary} if $\alpha$ vanishes on $\ao$,

\item \emph{real} if $\alpha$ vanishes on $\tc$,

\item \emph{complex} otherwise.\end{itemize}

Also, an imaginary root  $\alpha \in \Phi$ is called

\begin{itemize}
\item \emph{compact} if $\g_\alpha\subset \ka$, \emph{i.~e.} if 
$\theta$ is the $\id$ on $\g_\alpha$,

\item \emph{non-compact} if $\g_\alpha\subset \p$, \emph{i.~e.} if 
$\theta$ is $-\id$ on $\g_\alpha$.
\end{itemize}

A \emph{Vogan diagram} consists of the following data: a Dynkin diagram 
$\Delta$, an involution
$\dynkin \in \Aut (\Delta)$, and a subset $P$ of roots in 
$\Delta^{\dynkin}$.
The roots in $P$ are ``painted" in the graphical description of the 
Vogan diagram.
The Vogan diagram of a simple real Lie algebra $\go$ is as follows:
the Dynkin diagram corresponds to the complexification $\g$. Then one 
fixes a Cartan involution
$\theta_0$ and takes a maximally compact $\theta$-stable Cartan 
subalgebra $\ho$, together
with a system of simple roots $\Delta \subset \Phi(\g, \h)$
which is stable under the transpose of $\theta$;
$\dynkin$ is the restriction of ${}^t\theta$.  Finally, $P$ is the set 
of non-compact simple roots.

\bigbreak
A \emph{normalized Vogan diagram} is a Vogan diagram with at most one 
painted vertex,
\emph{i.~e.} at most one non-compact simple root.

\bigbreak The notion of normalized Vogan diagram helps to take care
of redundancies in the classification of real simple Lie algebras.
Indeed, the Theorem of Borel and de Siebenthal \cite[Th. 6.96]{Kn}
allows to go from the Vogan diagram of a simple real Lie algebra $\go$
to another one with at most one painted simple root, just by changing
appropriately the choice of the system of simple roots.

\bigbreak
Let now $\h$ be a fixed Cartan subalgebra of $\g$
and let $\Delta \subset \Phi(\g, \h)$ be a fixed system of simple 
roots.
Assume that $\invo$ is a sesquilinear Lie algebra involution of the 
form
$\tdynkin$ or $\jdynkin$, for $\dynkin: \Delta \to \Delta$ an
automorphism of the Dynkin diagram of order 1 or 2, and $J$ any
subset of the set $\Delta^{\dynkin}$.
Let $\go$ be the real form
$\go = \g^{\invo}$ of $\g$. Let $\ho := \h \cap \go$, a Cartan
subalgebra of $\go$. Our goal is to determine $\go$, and the 
isomorphism
classes of pairs $(\go, \ho)$.

Since $\invo$ commutes with $\omega$, $\omega$ preserves $\go$.
Let $\theta_0: \go \to \go$ be the linear Lie algebra involution
given by the restriction of $\omega$. Then $\theta_0$ is a Cartan
involution of $\go$ and $\theta_0(\ho)= \ho$; its complexification
is clearly $\theta = \invo\omega$. If $\sigma = \jdynkin$ the
transpose of $\theta$ preserves $\Delta$, and in fact coincides
with $\dynkin$; if $\sigma = \tdynkin$ the transpose of $\theta$
coincides with $-\dynkin$.

\begin{lema}\label{lematablados} Assume that $P = \Delta^{\dynkin}
- J$ has at most one element. Then the pair $(\go, \ho)$
determined by $\invo$ is as in Table \ref{tablados}.
\end{lema}

\begin{table}[t]
\begin{center}
\begin{tabular}{|p{1cm}|p{3,8cm}|p{3,4cm}|p{4,2cm}|p{3,7cm}|}
\hline {\bf Type} & {\bf Involution $\invo$} & $\go = \g^{\invo}$
& $\ho = \go \cap \h$
 & {\bf Remarks}
\\ \hline  All & $\tid$  & $\g_{\R}$ & $\sum_{\alpha\in \Delta} \R 
h_{\alpha}$  & Split form, split Cartan.
\\ \hline $A_{2n}$& $\tdynkin$, $\dynkin \neq  \id$   & $\mathfrak{su}(n, n  + 1)$
& $\sum_{\alpha\in \Delta^{\dynkin}} \R h_{\alpha}  +  \newline
\sum_{\alpha\in \Delta - \Delta^{\dynkin}} [\R (h_{\alpha} +
h_{\dynkin(\alpha)}) \newline + \R i(h_{\alpha} -
h_{\dynkin(\alpha)})]$ & $\ho$ is maximally non-compact.
\\ $A_{2n + 1}$& $\tdynkin$, $\dynkin \neq  \id$ &
$\go = \mathfrak{su}(n  + 1, n  + 1)$  & &
\\  $D_{n}$& $\tdynkin$, $\dynkin \neq  \id$   & $\mathfrak{so}(n  - 1, 
n+ 1)$  & &
\\  $E_{6}$& $\tdynkin$, $\dynkin \neq  \id$   & EII  & &
\\ \hline  All & $\comp$  & $\g_c$  & $\sum_{\alpha\in \Delta} i\R 
h_{\alpha}$ \newline & Compact form.
\\ \hline $A_{n}$& $\jid$, $\# P = 1$

& $\mathfrak{su}(j, n + 1 - j)$  & Idem.& $j\in P$
\\  $B_{n}$& $\jid$, $\# P = 1$
& $\mathfrak{so}(2j, 2n + 1 - 2j)$  & & $j\in P$
\\  $C_{n}$& $\jid$, $\# P = 1$
& $\mathfrak{sp}(j, n  - j)$ & & $j\in P$
\\  $D_{n}$& $\jid$, $\# P = 1$ \newline ${}^{}$\quad ${}^{}$\, $\# P = 
1$
& $\mathfrak{so}(2j, 2n - 2j)$
\newline $\mathfrak{so}^*(2n)$ & & $j\in P$, $j \le n- 2$\newline 
$j\in P$, $j > n- 2$
\\  $E_{6}$& $\jid$, $\# P = 1$
\newline ${}^{}$\newline ${}^{}$
& EII,  EIII  & & $P = \{$extrem of the short, respect. long
branch$\}$.
\\  $E_{7}$& $\jid$, $\# P = 1$
& EV, EVI, EVII   & & $P = \{$extrem of the short, resp. medium,
long,  branch$\}$.
\\  $E_{8}$& $\jid$, $\# P = 1$
& EVIII, EIX   & & $P = \{$extrem of the medium, resp. long,
branch$\}$.
\\ $F_{4}$& $\jid$ , $\# P = 1$
& FI, FII & & $P = \{$extrem,  short, respect. long$\}$.
\\  $G_{2}$& $\jid$, $\# P = 1$   & G   & & Split form.
\\ \hline $A_{2n}$& $\jdynkin$, $\dynkin \neq  \id$, $\# P = 0$   & $\mathfrak{sl}(2n+1, \R)$
& $\sum_{\alpha\in \Delta^{\dynkin}} \R i h_{\alpha}  +  \newline
\sum_{\alpha\in \Delta - \Delta^{\dynkin}} [\R i(h_{\alpha}  +
h_{\dynkin(\alpha)}) \newline + \R (h_{\alpha} -
h_{\dynkin(\alpha)})]$ &
\\  $A_{2n + 1}$& $\jdynkin$, $\dynkin \neq  \id$, $\# P = 0$
\newline ${}^{}$ \qquad ${}^{}$ \qquad ${}^{}$ \quad $\# P = 1$
& $\mathfrak{sl}(\frac{n+1}2, \mathbb H)$
\newline$\mathfrak{sl}(2n+2, \R)$ & &
\\  $D_{n}$& $\jdynkin$, $\dynkin \neq  \id$, $\# P = 0$
\newline ${}^{}$ \qquad ${}^{}$ \qquad ${}^{}$ \quad $\# P = 1$
&$\mathfrak{so}(1, 2n  - 1)$\newline  $\mathfrak{so}(2j +1, 2(n -
j)  - 1)$& ${}^{}$
\newline &$j\in P$
\\  $E_{6}$& $\jdynkin$, $\dynkin \neq  \id$, $\# P = 0$
\newline ${}^{}$ \qquad ${}^{}$ \qquad ${}^{}$ \quad $\# P = 1$
& EIV \newline EI  & &
\\ \hline
\end{tabular}
\end{center}

\

\caption{Determination of $(\go, \ho)$. See \cite{Kn} for  conventions 
on Vogan diagrams.} \label{tablados}
\end{table}

\pf Let us analyze the different possibilities for  $\invo$.

\bigbreak (i). $\invo = \tid$. Then $\ho$ is the real span of
$h_{\alpha}$, $\alpha \in \Delta$; it is a split Cartan subalgebra
of $\go$ and this last is the split real form of $\g$. Also,
$\dc(\ho) = 0$, $\dnc(\ho) = \# \Delta$.

\bigbreak (ii). $\invo = \tdynkin$, $\dynkin \neq \id$. Then $\ho
= \tc \oplus \ao$, where $\tc$ is the real span of $\{i(h_{\alpha}
- h_{\dynkin(\alpha)}): \alpha \in \Delta - \Delta^{\dynkin}\}$;
and $\ao$ is the real span of $\{h_{\alpha}: \alpha \in
\Delta^{\dynkin}\} \cup$ $\{h_{\alpha} + h_{\dynkin(\alpha)}:
\alpha \in \Delta - \Delta^{\dynkin}\}$. 
In this case, $\invo$ is a \emph{Steinberg normal form, i.~e.}, it
leaves a Borel subalgebra invariant. The list of all such is well-known, see for
example \cite{W1, W2}. For the sake of completeness we briefly complete the argument.

It is easy to see that
there are no imaginary roots; therefore, $\ho$ is maximally
non-compact by \cite[Prop. 6.70]{Kn}. Also,  $\dc(\ho) = \dfrac
{\#(\Delta - \Delta^{\dynkin})}2$, $\dnc(\ho) = \# \Delta - \dfrac
{\#(\Delta - \Delta^{\dynkin})}2 = \# \Delta^{\dynkin} + \dfrac
{\#(\Delta - \Delta^{\dynkin})}2 =$ rank of the associated
symmetric space. We can not apply the results in \cite[Chapter
VI]{Kn}. However, the necessary information can be obtained
looking at the classification of symmetric spaces in \cite[Table V, p. 
518]{H}.

\bigbreak If $\g$ is of type $A_{2n}$, then the rank of the
associated symmetric space is $n$. We conclude that $\go =
\mathfrak{su}(n, n  + 1)$. If $\g$ is of type $A_{2n + 1}$, then
the rank of the associated symmetric space is $n + 1$. We conclude
that $\go = \mathfrak{su}(n + 1, n  + 1)$.  
If $\g$ is of type $D_n$, $n \ge 4$, then the rank of
the associated symmetric space is $n-1$. We conclude that $\go =
\mathfrak{so}(n - 1, n  + 1)$. If $\g$ is of type $E_{6}$, then
the rank of the associated symmetric space is $4$. We see that
$\go$ is  of type EII.

\bigbreak (iii). $\invo = \jid$. Then $\ho$ is the real span of $i
h_{\alpha}$, $\alpha \in \Delta$; it is clearly maximally compact
and $\Delta$ is $\theta$-stable. It is easy to see that  $P =
\Delta - J$, \emph{i.~e.} that the compact roots are precisely
those in $J$. Also,  $\dnc(\ho) = 0$, $\dc(\ho) = \# \Delta$.
 We can then apply the results in \cite[Chapter VI, pages 355 to 
362]{Kn}.

\bigbreak (iv). $\invo = \jdynkin$, $\dynkin \neq \id$. Then $\ho$
is the real span of $\{i h_{\alpha}: \alpha \in \Delta^{\dynkin}\}
\cup$ $\{i(h_{\alpha} + h_{\dynkin(\alpha)}): \alpha \in \Delta -
\Delta^{\dynkin}\} \cup$ $\{h_{\alpha} - h_{\dynkin(\alpha)}:
\alpha \in \Delta - \Delta^{\dynkin}\}$. Clearly, $\Delta$ is
$\theta$-stable. It is easy to see that there are no real roots;
therefore, $\ho$ is maximally compact by \cite[Prop. 6.70]{Kn}. It
is easy to see that  $P = \Delta^{\mu} - J$, \emph{i.~e.} that
the compact roots are precisely those in $J$. Also,  $\dnc(\ho) =
\dfrac {\#(\Delta - \Delta^{\dynkin})}2$, $\dc(\ho) = \# \Delta -
\dfrac {\#(\Delta - \Delta^{\dynkin})}2$. We can then apply the
results in \cite[Chapter VI, pages 355 to 362]{Kn}.
 \epf

\subsection{Lie bialgebras}
Recall that a quasitriangular Lie bialgebra $(\g, r)$ is called
\emph{factorizable}  if $r + r^{21} \in S^2 \g$
defines a nondegenerate inner product on $\g^*$ \cite{RS}.

\begin{definition}\label{alm-fact}
We shall say that a real Lie bialgebra $(\lgo, \delta)$ is \emph{almost 
factorizable} if the complexification $(\lgot, \delta)$ is factorizable.

\bigbreak 
We shall consider the following  particular class of almost factorizable Lie bialgebras. 
A  Lie bialgebra $(\lgo, \delta)$ is \emph{imaginary factorizable}
if the complexification $(\lgot, \delta)$ is factorizable and 
$r\in \lgot \otimes \lgot$ is given by
\begin{equation}\label{almost-factorizable}
r = \rla + i \rom, \quad \text{where }\rla \in \Lambda^2(\lgo), \quad \rom \in S^2(\lgo). 
\end{equation}
\end{definition}

Real factorizable Lie bialgebras are of course almost factorizable, but 
we shall see that the converse is not true. In fact, real simple Lie bialgebras
are triangular, factorizable or imaginary factorizable.

\bigbreak 
Doubles and duals of imaginary factorizable Lie bialgebras are 
computed in Proposition \ref{im-double} below.

\begin{lema}\label{whead} Let $(\go, \delta)$ be an absolutely simple 
real Lie bialgebra.
Then $\delta  = \partial \ro$ for a unique $\ro \in \Lambda^2(\go)$ and
$(\go, \delta)$ is either triangular or almost factorizable.
\end{lema}
\pf Let $\g$ be the complexification of $\go$.
By Whitehead's Lemma, there exists $\ro \in \Lambda^2(\go)$ such that 
$\delta(x) = \ad _x \ro$
for all $x\in \go$, and \emph{a fortiori} for all $x\in \g$.
(Note that $\ro$ is uniquely defined by $\delta$ since 
$\Lambda^2(\go)^{\go} = 0$).
Then $[\ro^{12}, \ro^{13}] + [\ro^{12}, \ro^{23}] + [\ro^{13}, 
\ro^{23}] \in \Lambda^3(\go)^{\go}$,
and there exists $c\in \Cc$ such that
$[\ro^{12}, \ro^{13}] + [\ro^{12}, \ro^{23}] + [\ro^{13}, \ro^{23}]= 
c^2 [\Omega^{13}, \Omega^{23}]$.
Since $\Omega \in \go \otimes \go$,
$c^2 \in \R$. If $c = 0$, $(\go, \delta)$ is triangular. Assume $c \neq 
0$, let $t = 2ic$ and let
\begin{equation}\label{r}
r = \ro + ic\Omega = \ro + \dfrac{t}2\Omega.
\end{equation}
Then $(\g, \delta)$ is quasitriangular, with $\delta = \partial r$;
and it is furthermore factorizable since $r + r^{21} = 2ic\Omega$.
\epf

Let $t$ be as in the proof of the Lemma; it is  defined up to a sign, 
since $c^2$ is unique.

\bigbreak 
If $t\in \R - 0$, then $(\go, \delta)$ is quasitriangular with $r$ 
given by \eqref{r}.

\bigbreak 
If $t\in i\R - 0$, then $(\go, \delta)$ is \emph{not} quasitriangular,
but \emph{imaginary factorizable}.

\bigbreak 
The distinction between these cases  is already present in \cite{CGR, Ch, 
P1, S}, from a different point
of view.
There are quantizations of these Lie bialgebras: these are
$*$-Hopf algebras in the imaginary factorizable case, and real
Hopf algebras in the factorizable case.

\subsection{The theorem of Belavin and Drinfeld}\label{BDsection}
Let $\g$ be a complex simple Lie algebra and
let $\h \subset \g$ be a Cartan subalgebra. Let $\Delta$ be a choice of
a set of simple roots.

\begin{definition}
A Belavin-Drinfeld triple (BD-triple for short) is a triple
$(\Gamma_1, \Gamma_2, T)$ where $\Gamma_1,$ $\Gamma_2$ are subsets of 
$\Delta$ and
$T: \Gamma_1 \to \Gamma_2$ is a bijection that preserves the inner 
product
and satisfies the nilpotency condition:
for any $\alpha \in \Gamma_1$ there exists a positive integer $n$ for 
which
$T^n (\alpha)$ belongs to $\Gamma_2$ but not to $\Gamma_1$.
\end{definition}
Let $(\Gamma_1, \Gamma_2, T)$  be a Belavin--Drinfeld triple.
Let $\widehat{\Gamma_i}$ be the set of positive roots
lying in the subgroup generated by $\Gamma_i$, for $i = 1,2$.
There is an associated partial ordering on $\Phi^+$ given by
$\alpha \prec \beta$ if
$\alpha \in \widehat{\Gamma_1}$, 
$\beta \in \widehat{\Gamma_2}$, and $\beta = T^n(\alpha)$
for a positive integer $n$.

A \emph{continuous parameter} for the Belavin--Drinfeld triple 
$(\Gamma_1, \Gamma_2, T)$
is an element $\lambda \in \h^{\otimes 2}$ such that
\begin{align}
(T (\alpha) \otimes 1) \lambda + (1 \otimes \alpha) \lambda &= 0,
\quad \text{for all } \alpha \in \Gamma_1,
 \\
\lambda + \lambda^{21} &= \Omega_0.
\end{align}

Let $\ag_1$, $\ag_2$ be the reductive subalgebras of $\g$ with Cartan 
subalgebras
generated by $h_{\alpha}$, $\alpha$ in $\Gamma_1$, resp. in $\Gamma_2$, 
and
with Dynkin diagrams $\Gamma_1$, resp. in $\Gamma_2$. We extend $T$ to
a Lie algebra isomorphism $\widehat T: \ag_1\to \ag_2$.

\begin{theorem}{\em{(}}Belavin--Drinfeld, see \cite{BD}{\em{\/)}}.
Let $(\g, \delta)$ be  a factorizable complex simple
Lie bialgebra. Then there exist a Cartan subalgebra $\h$,
a system of simple roots $\Delta$, a Belavin--Drinfeld
triple $(\Gamma_1, \Gamma_2, T)$,  a continuous parameter $\lambda$
and $t\in \Cc - 0$ such that the $r$-matrix is given by
\begin{equation}\label{BD}
r =t\left(\lambda + \sum_{\alpha \in \Phi^+} x_{-\alpha} 
\otimes x_\alpha
        + \sum_{\alpha, \beta \in \Phi^+,  \alpha \prec \beta}
                               x_{-\alpha} \wedge x_\beta\right),
\end{equation}
where $x_\beta \in \g_\beta$, $\beta \in \pm\Phi^{+}$, are root vectors  
normalized by
\begin{align}\label{norm1}
(x_\beta \vert x_{-\beta}) = 1, \qquad \text{for all }\beta \in \Phi^+,
\\ \label{norm2}
\widehat T(x_\beta) = x_{T(\beta)}, \qquad \text{for all }\beta \in 
\Gamma_1.
\end{align}
\qed
\end{theorem}
Clearly, $r + r^{21} = t\Omega$.

\section{Proof of the main theorem}

The proof of the main result is split into two lemmas.

\begin{lema}\label{necesario}

(I). Let $(\go, \delta)$ be an absolutely simple real Lie bialgebra.
Let $\g$ be the complexification of $\go$ and let $\invo$ be the 
sesquilinear involution
of $\g$ whose fixed-point set is $\go$. Assume that $(\go, \delta)$ is 
almost factorizable.
Then there exist:

\begin{itemize}
\item a complex number $c \neq 0$ with $c^2 \in \R$; set $t = 2ic$ and fix the 
choice $t\in \R_{>0} \cup i\R_{>0}$;
\item a Cartan subalgebra $\h$ of $\g$;
\item a system of simple roots $\Delta\subset \Phi(\g, \h)$;
\item a Belavin-Drinfeld triple $(\Gamma_1, \Gamma_2, T)$
and a continuous parameter $\lambda \in \h^{\otimes 2}$;
\end{itemize}
such that \begin{itemize}
\item $\h$ is stable under $\invo$ (we denote $\ho := \h \cap \go$).
\item $\invo^*(\Delta)$ is either $\Delta$ or $-\Delta$;
furthermore $\dynkin: = \invo^*: \Delta \to \pm \Delta$ is an
automorphism of the Dynkin diagram. If $\invo^*(\Delta) = \Delta$
then there exists an automorphism $\dynkin$ of the Dynkin diagram
such that, for an appropriate choice of the $e_{\alpha}$'s,
$\invo$ is $\tdynkin$.

If $\invo^*(\Delta) = -\Delta$ then there exists an automorphism
$\dynkin$ of the Dynkin diagram and $J \subset \Delta^{\dynkin}$
such that, for an appropriate choice of the $e_{\alpha}$'s,
$\invo$ is  $\jdynkin$.

\item $\delta = \partial \ro$ where $\ro \in \Lambda^2(\go)$ is given 
by the formula
\begin{equation}\label{BDreal}
\ro =t\left(\frac 12(\lambda - \lambda^{21}) + \frac 12
\sum_{\alpha \in \Phi^+} x_{-\alpha} \wedge x_\alpha
        + \sum_{\alpha, \beta \in \Phi^+,  \alpha \prec \beta}
                               x_{-\alpha} \wedge x_\beta\right).
\end{equation}\end{itemize}

\bigbreak (II). Let $(\go, \delta)$, $(\go', \delta')$ be two
almost factorizable absolutely simple real Lie bialgebras. Let
$\psi: \go' \to \go$ be an isomorphism of Lie bialgebras. Let
$\invo$, $\invo'$ be the involutions corresponding to $\go$,
$\go'$. Let $\h'$, $\ho'$, $\Delta'\subset \Phi(\g, \h')$;
$(\Gamma_1', \Gamma_2', T')$, $\lambda' \in {\h'}^{\otimes 2}$,
$c'$, etc., be the corresponding objects for the Lie bialgebra
$\go'$. Then
\begin{itemize}
\item $\psi(\h') = \h$.
\item If $\invo^*(\Delta) = \Delta$, resp.  $-\Delta$, then
${\invo'}^*(\Delta') = \Delta'$, resp.  $-\Delta'$.
\item $\psi$ induces an isomorphism of Dynkin diagrams $\psi^*: \Delta
\to \Delta'$;
\item $\mu'\psi^* = \psi^*\mu$, $J' = \psi^*(J)$;
\item $(\Gamma'_1, \Gamma'_2, T') = (\psi^*(\Gamma_1), 
\psi^*(\Gamma_2),
\psi^*T{\psi^*}^{-1})$.
\end{itemize}
\end{lema}

In other words, the Lemma says that, for an appropriate choice of
the $e_{\alpha}$'s, $\invo$ is either $\tid$ or $\tdynkin$ with
$\mu\neq \id$, or $\jid$, or $\comp$, or $\jdynkin$ with $\mu\neq
\id$; and that this does not depend on the isomorphism class as
Lie bialgebra.

The uniqueness in Part II uses in an essential way that the
parameter $t$ is in  $\R_{>0} \cup i\R_{>0}$. Changing $t$ to $-t$
would affect \eqref{BDreal} by changing $\Delta$ to $-\Delta$,
etc.

\pf (I). By the Theorem of Belavin and Drinfeld, there exist a
Cartan subalgebra $\h$ of $\g$, a system of simple roots $\Delta$,
a Belavin-Drinfeld triple $(\Gamma_1, \Gamma_2, T)$, a continuous
parameter $\lambda \in \h^{\otimes 2}$ and a non-zero complex number $c$
(with $c^2 \in \R$ by Lemma \ref{whead}) such that the
complexification $(\g, \delta)$ is quasitriangular with $r$-matrix
$r$ given by \eqref{BD}  and $t = 2ic$. Then $\delta = \partial
\ro$, where $\ro$ is given by \eqref{BDreal}; and $\ro \in
\Lambda^2(\go)$ by the uniqueness in Lemma \ref{whead}.

Let $H$ be the image of $\ro$ under the Lie bracket $[\,,\,]: \go
\otimes \go \to \go$; then $H = -2ic \sum_{\alpha \in \Phi^+}
h_{\alpha}$. Indeed, if $\alpha \prec \beta$ then $[x_{-\alpha} ,
x_\beta]  = 0$ because $\alpha \neq \beta$ but both have the same
level with respect to $\Delta$. Note that here the nilpotency
condition on the BD-triple is used. Since $H$ is a regular element
of $\h$ by \cite[13.3]{Hu}, $\h = \cent_{\g} (H)$ is a Cartan
subalgebra of $\g$; since $\h$  is the complexification of
$\cent_{\go} (H)$, we see that $\h$ is stable under $\invo$. By
Lemma \ref{involuciones},  $\Phi(\g, \h)$ is stable under
$\invo^*$, and $\invo(\g_{\alpha}) = \g_{\invo^*(\alpha)}$, for all
$\alpha \in \Phi(\g, \h)$.

Let $$W = \sum_{\alpha, \beta \in \Phi: \alpha + \beta = 0} \g_{\alpha} 
\otimes \g_{\beta}, \qquad
U = \sum_{\alpha, \beta \in \Phi: \alpha + \beta \neq 0} \g_{\alpha} 
\otimes \g_{\beta}.$$
It is clear that
$$
(\invo\otimes\invo)\left(\sum_{\alpha \in \Phi^+} x_{-\alpha} \wedge 
x_\alpha\right) \in W, \qquad
(\invo\otimes\invo)\left(\sum_{\alpha, \beta \in \Phi^+,  \alpha \prec 
\beta}
 x_{-\alpha} \wedge x_\beta\right) \in U.
$$
Assume that $c \in i\R$, \emph{i.~e.} that $(\go, \delta)$ is 
quasitriangular.
It follows that $(\invo\otimes\invo)(\sum_{\alpha \in \Phi^+} 
x_{-\alpha} \wedge x_\alpha)
= \sum_{\alpha \in \Phi^+} x_{-\alpha} \wedge x_\alpha$. Since the 
elements
$x_{-\alpha} \wedge x_\alpha$ are linearly independent in $W$, we 
conclude that
$\invo^*(\Delta) = \Delta$; $\dynkin$ is an automorphism of the Dynkin 
diagram
by Lemma \ref{involuciones}.

Similarly, if $c \in \R$, \emph{i.~e.} if $(\go, \delta)$ is imaginary
factorizable, we conclude that $\invo^*(\Delta) = -\Delta$ and
$\dynkin$ is an automorphism of the Dynkin diagram
again by Lemma \ref{involuciones}.

\bigbreak (II). We assume that $\g$ is the complexification of
both $\go$, $\go'$ (equality, not just isomorphism); $\psi_0$
extends to a Lie algebra automorphism $\psi$ of $\g$, and
$\invo\psi = \psi\invo'$.

Let $\ro \in \Lambda^2(\go)$, $\ro' \in \Lambda^2(\go')$ be such
that $\delta =
\partial \ro$, $\delta' = \partial \ro'$.  Thus
\begin{equation*}
\ro' = t' \left( \frac 12(\lambda' - {\lambda'}^{21}) +  \frac 12
\sum_{\alpha \in \Phi(\g, \h')^+} x'_{-\alpha} \wedge x'_\alpha
        + \sum_{\alpha, \beta \in \Phi(\g, \h')^+,  \alpha \prec' 
\beta}
                               x'_{-\alpha} \wedge x'_\beta\right).
\end{equation*}

Since $\ro$ is unique, $(\psi_0 \otimes \psi_0) (\ro') =  \ro$ and
\emph{a fortiori} $\psi (\ho') =  \ho$, $\psi(\h') =  \h$. Thus
\begin{equation*}
\ro =t' \left( \frac 12 (\psi \otimes \psi) ( \lambda' -
{\lambda'}^{21} ) + \frac 12 \sum_{\alpha' \in \Phi(\g, \h')^+}
\psi(x'_{-\alpha'}) \wedge \psi(x'_{\alpha'})
        + \sum_{\alpha', \beta' \in \Phi(\g, \h')^+,  \alpha' \prec' 
\beta'}
                               \psi(x'_{-\alpha'}) \wedge 
\psi(x'_{\beta'})\right).
\end{equation*}

But also $\psi$ induces a bijection $\psi^*: \Phi(\g, \h) \to
\Phi(\g, \h')$, and $\psi(\g_{\psi^*\alpha}) = \g_{\alpha}$.
Arguing as above, we conclude that

\begin{align}
\label{aweauno}t' (\psi \otimes \psi) (\lambda' -
{\lambda'}^{21} ) &=t(\lambda - \lambda^{21}),\\
\label{aweados}t'\sum_{\alpha' \in \Phi(\g, \h')^+}
\psi(x'_{-\alpha'}) \wedge \psi(x'_{\alpha'}) &= t
\sum_{\alpha \in \Phi^+} x_{-\alpha} \wedge x_\alpha,\\
\label{aweatres}t'\sum_{\alpha', \beta' \in \Phi(\g,
\h')^+, \alpha' \prec' \beta'} \psi(x'_{-\alpha'}) \wedge
\psi(x'_{\beta'}) &=t\sum_{\alpha, \beta \in \Phi^+,
\alpha \prec \beta}
                               x_{-\alpha} \wedge x_\beta.
\end{align}

Now, $(\go, \delta)$, $(\go', \delta')$ are either both
quasitriangular, or both imaginary factorizable; that is, are
either both $t$, $t'$ are in $\R_{> 0}$, or both in $i\R_{> 0}$.
By \eqref{aweados}, $\Phi(\g, \h')^+ = \psi^*(\Phi(\g, \h)^+)$ and
hence $\Delta' = \psi^*(\Delta)$.

Hence, \eqref{aweatres} implies that $(\Gamma'_1, \Gamma'_2, T') =
(\psi^*(\Gamma_1), \psi^*(\Gamma_2), \psi^*T{\psi^*}^{-1})$.

From the equality ${\invo'}^*\psi^* = \psi^*\invo^*$ we conclude
that $\mu'\psi^* = \psi^*\mu$.

Finally, recall from the proof of the Lemma \ref{involuciones}
that $J = \{\alpha \in \Delta^{\mu}: \sigma(e_{\alpha}) =
c_{\alpha} e_{\invo^*\alpha}$ with $c_{\alpha} < 0\}$. It follows
without difficulties that $J' = \psi^*(J)$. \epf

\begin{obs} The description of $\ro$ in \eqref{BDreal} depends on
the choice of a family $x_{\alpha} \in \g_{\alpha}$, $\alpha \in \pm 
\Phi^{+}$,
satisfying \eqref{norm1}, \eqref{norm2}. Such a family can be 
constructed
starting from any choice of $x_{\alpha} \in \g_{\alpha}$ for $\alpha 
\in \Gamma_1 - \Gamma_2$;
different choices do not affect \eqref{BDreal}.

On the other hand, the choice of elements $e_{\alpha}$ in Lemma 
\ref{involuciones} is independent of the $x_{\beta}$'s; the arguments in the preceding and 
next lemmas do not depend on the explicit form of the $x_{\alpha}$ but on which 
root space they are living in.
\end{obs}

\begin{definition}\label{dynkinstable}  Let $\dynkin$ be an 
automorphism of the Dynkin diagram.
A BD-triple  $(\Gamma_1, \Gamma_2, T)$ is \emph{$\dynkin$-stable}
if $\dynkin(\Gamma_1) =\Gamma_1$, $\dynkin(\Gamma_2) =\Gamma_2$,
and $T\dynkin = \dynkin T$. A BD-triple  $(\Gamma_1, \Gamma_2, T)$
is \emph{$\dynkin$-antistable} if $\dynkin(\Gamma_1) =\Gamma_2$,
$\dynkin(\Gamma_2) =\Gamma_1$, and $T^{-1}\dynkin = \dynkin T$. In
particular, if $\dynkin = \id$ then all BD-triples are
$\dynkin$-stable, and the only BD-triple $\dynkin$-antistable has
$\Gamma_1 =  \Gamma_2 = \emptyset$.
\end{definition}

\begin{lema}\label{varsigma}
Let $\g$ be a simple complex Lie algebra, $\h$ a Cartan subalgebra  of 
$\g$,
$\Delta\subset \Phi(\g, \h)$ a system of simple roots,
$(\Gamma_1, \Gamma_2, T)$  a BD-triple,
$\lambda \in \h^{\otimes 2}$ a continuous parameter and $t$ a complex 
number.
Write $$\lambda - \lambda^{21} = 
\sum_{\alpha, \beta \in \Delta} \lambda_{\alpha, \beta} h_{\alpha} 
\wedge h_\beta.$$
By convention, $\lambda_{\alpha, \beta} = -\lambda_{\beta, \alpha}$ for all $\alpha, \beta \in \Delta$.
Let $\ro$ be given by formula \eqref{BDreal}. Let $\invo$ be an 
involution of the form
$\tid, \comp, \jid$, or $\tdynkin$, $\jdynkin$ with $\dynkin \neq \id$.
Let $\go$ be the real Lie algebra of vectors fixed by $\invo$.
Then:

\medbreak (a). Assume that $\invo = \tid$. Then $\ro \in
\Lambda^2(\go)$ if and only if $t\in \R$, $\lambda_{\alpha,
\beta}\in \R$ for all $\alpha, \beta \in \Delta$ (no restrictions
on the BD-triple).

\medbreak (b).  Assume that $\invo = \tdynkin$. Then $\ro \in 
\Lambda^2(\go)$ if and only if
$t\in \R$, $\lambda_{\alpha, \beta} = 
\overline{\lambda_{\dynkin(\alpha), \dynkin(\beta})}$
for all $\alpha, \beta \in \Delta$ and
the BD-triple is $\dynkin$-stable.

\medbreak (c).  Assume that $\invo = \jid$, or $\invo = \comp$. Then 
$\ro \in \Lambda^2(\go)$ if and only if
$t\in i\R$, $\lambda_{\alpha, \beta}\in i\R$, for all $\alpha, \beta 
\in \Delta$ and
the BD-triple has $\Gamma_1 =  \Gamma_2 = \emptyset$.

\medbreak (d).  Assume that $\invo = \jdynkin$. Then $\ro \in 
\Lambda^2(\go)$ if and only if
$t\in i\R$, $\lambda_{\alpha, \beta} = 
-\overline{\lambda_{\dynkin(\alpha), \dynkin(\beta})}$,
for all $\alpha, \beta \in \Delta$ and
the BD-triple is $\dynkin$-antistable.

\end{lema}

\pf Since $\invo \otimes \invo$ is the sesquilinear involution of
$\g \otimes \g$ corresponding to $\go\otimes \go$, it is clear
that $\ro \in \Lambda^2(\go)$ if and only if $(\invo \otimes
\invo) (\ro) = \ro$. The lemma follows by a direct computation
that we include below for completeness. Note that $$(\invo \otimes
\invo) (\ro) = \overline{t} \left( \frac 12 (\invo \otimes
\invo)(\lambda - \lambda^{21}) + \frac 12 \sum_{\alpha \in \Phi^+}
(\invo \otimes \invo) (x_{-\alpha} \wedge x_\alpha) +
\sum_{\alpha, \beta \in \Phi^+,  \alpha \prec \beta} (\invo
\otimes \invo) (x_{-\alpha} \wedge x_\beta)\right) = \ro$$ if and
only if
\begin{align}
\label{conduno}
\overline{t} \sum_{\alpha, \beta \in \Delta} \overline{\lambda_{\alpha, 
\beta}}
\invo (h_{\alpha}) \wedge\invo (h_\beta)
&= t   \sum_{\alpha, \beta \in \Delta} \lambda_{\alpha, \beta} 
h_{\alpha} \wedge h_\beta,
\\ \label{conddos}
\overline{t} \sum_{\alpha \in \Phi^+} \invo (x_{-\alpha}) \wedge \invo 
(x_\alpha)
&= t  \sum_{\alpha \in \Phi^+} x_{-\alpha} \wedge x_\alpha,
\\ \label{condtres}
\overline{t}\sum_{\alpha, \beta \in \Phi^+,  \alpha \prec \beta} \invo 
(x_{-\alpha}) \wedge \invo (x_\beta)
&= t \sum_{\alpha, \beta \in \Phi^+,  \alpha \prec \beta} x_{-\alpha} 
\wedge x_\beta,
\end{align}
as in the proof of Lemma \ref{necesario}. In cases (a) and (b),
\eqref{conddos} holds if and only if $t\in \R$. Then
\eqref{conduno} holds if and only if $\lambda_{\alpha, \beta}\in
\R$ --in case (a)-- or $\lambda_{\alpha, \beta} =
\overline{\lambda_{\dynkin(\alpha), \dynkin(\beta})}$ --in case
(b)-- for all $\alpha, \beta \in \Delta$. Finally, in case (b),
\eqref{condtres} holds if and only if $\sum_{\alpha, \beta \in
\Phi^+,  \alpha \prec \beta} x_{-\dynkin(\alpha)} \wedge
x_{\dynkin(\beta)} =  \sum_{\alpha, \beta \in \Phi^+,  \alpha
\prec \beta} x_{-\alpha} \wedge x_\beta$. 

\bigbreak It is then easy to see that equality \eqref{condtres} holds 
if and
only if the BD-triple is $\dynkin$-stable. Indeed, let $\mathcal P =
\{(\alpha, \beta) \in \Phi^+ \times \Phi^+:  \alpha \prec
\beta\}$. The equality implies that $\dynkin \times \dynkin (\mathcal P) =
\mathcal P$ and that $\dynkin(\alpha) \prec \dynkin(\beta)$ if $\alpha
\prec \beta$. Since $\alpha \prec T\alpha$ for any $\alpha \in
\Gamma_1$, we conclude that $\dynkin(\Gamma_1) =\Gamma_1$, and
similarly, $\dynkin(\Gamma_2) =\Gamma_2$. For any $\alpha \in
\Gamma_1$, we set $$ X(\alpha) = \{\gamma \in \Gamma_1: \alpha
\prec \gamma \text{ or } \gamma \prec \alpha \text{ or } \gamma =
\alpha\} = \{\gamma \in \Gamma_1:  \gamma = T^m\alpha, \quad m\in
\Z\}.$$ Clearly, $\dynkin(X(\alpha)) = X(\dynkin(\alpha))$. 
If $\card X(\alpha) = 1$, it is easy to see that $T\dynkin(\alpha) 
= \dynkin T(\alpha)$. Assume that $\card X(\alpha) > 1$. Now,
there exists $\alpha_0$ such that $X(\alpha) = \{\alpha_0,
\alpha_1, \dots, \alpha_s\}$, where $\alpha_i \prec \alpha_{i+1} =
T\alpha_i$, $0\le i < s$. Thus, $\dynkin(X(\alpha)) =
\{\dynkin(\alpha_0), \dynkin(\alpha_1), \dots,
\dynkin(\alpha_s)\}$ and $T\dynkin(\alpha_i) =
\dynkin(\alpha_{i+1}) = \dynkin(T\alpha_i)$ for all $i$. This
shows that $T\dynkin = \dynkin T$.

\bigbreak
In cases (c) and (d), \eqref{conddos} holds if
and only if $t\in i\R$. Then \eqref{conduno} holds if and only if
$\lambda_{\alpha, \beta}\in i\R$ --in case (c)-- and
$\lambda_{\alpha, \beta} = -\overline{\lambda_{\dynkin(\alpha),
\dynkin(\beta})}$ --in case (d)-- for all $\alpha, \beta \in
\Delta$. Finally, in case (d), \eqref{condtres} holds if and only
if $\sum_{\alpha, \beta \in \Phi^+,  \alpha \prec \beta}
x_{\dynkin(\alpha)} \wedge x_{-\dynkin(\beta)} =  \sum_{\alpha,
\beta \in \Phi^+,  \alpha \prec \beta} x_{-\alpha} \wedge
x_\beta$.

\bigbreak It is easy then to see that equality \eqref{condtres} holds 
if and
only if the BD-triple is $\dynkin$-antistable. Indeed, let $\mathcal P^t =
\{(\beta, \alpha) \in \Phi^+ \times \Phi^+:  \alpha \prec
\beta\}$. The equality implies that $\dynkin \times \dynkin (\mathcal P) =
\mathcal P^t$ and that $\dynkin(\beta) \prec \dynkin(\alpha)$ if $\alpha
\prec \beta$. Since $\alpha \prec T\alpha$ for any $\alpha \in
\Gamma_1$, we conclude that $\dynkin(\Gamma_1) =\Gamma_2$, and
similarly, $\dynkin(\Gamma_2) =\Gamma_1$. Keep the notation for
$X(\alpha)$ as above. For any $\beta \in \Gamma_2$, we set $$
Y(\beta) = \{\delta \in \Gamma_2: \beta \prec \delta \text{ or }
\delta \prec \beta \text{ or } \delta = \beta\} = \{\delta \in
\Gamma_2:  \delta = T^m\beta, \quad m\in \Z\}.$$ Clearly,
$\dynkin(X(\alpha)) = Y(\dynkin(\alpha))$. If $\card X(\alpha) = 1$, 
it is easy to see that $T^{-1}\dynkin(\alpha) 
= \dynkin T(\alpha)$. Assume that $\card X(\alpha) > 1$. Now, there 
exists
$\beta_0$ such that $Y(\dynkin(\alpha)) = \{\beta_0, \beta_1,
\dots, \beta_s\}$, where $\beta_i \prec \beta_{i+1} = T\beta_i$,
$0\le i < s$. Thus, $\dynkin(\alpha_i) = \beta_{s - i}$ and
$\dynkin(T\alpha_i) = \dynkin(\alpha_{i+1})= \beta_{s - i -1} =
T^{-1}\beta_{s - i} = T^{-1}\dynkin(\alpha_i)$ for all $i$. This
shows that $T^{-1}\dynkin = \dynkin T$. \epf

\bigbreak
We have collected now all the necessary information to prove the
main result.

\bigbreak
 \emph{Proof of Theorem \ref{main1}. \/} By Lemmas
\ref{necesario} part (I) and \ref{varsigma}, we know the existence
of $\h$, $\Delta$, $\invo$ of the type $\tdynkin$ or $\jdynkin$,
$(\Gamma_1, \Gamma_2, T)$, $t\in \R_{>0} \cup i\R_{>0}$, and
$\lambda \in \h^{\otimes 2}$; such that $(\go, \delta)$ is
isomorphic as real Lie bialgebra to $(\g^{\invo},
\partial \ro)$, where $\ro$ is given by \eqref{BDreal}. 

It remains to determine when different data in Table \ref{tablauno} give rise
to isomorphic Lie bialgebras. By Lemma \ref{necesario} part (II) and Lemma \ref{lematablados}, 
isomorphic Lie bialgebras can arise only from data in the same row; then the statement follows
from Lemma \ref{necesario} part (II). \qed

\begin{obs}\label{isom} Let $\ho = \cent_{\go} (H)$, where $H$ is the image of
$r_0$ under the Lie bracket as in the proof of Lemma \ref{necesario}.
If $\invo = \jdynkin$, or if $\invo = \jid$, then $\ho$ is a maximally compact
Cartan subalgebra of $\go$ and Lemma \ref{necesario} provides a Vogan diagram; however,
this Vogan diagram is not normalized. Even if  the Theorem of Borel and de Siebenthal \cite[Th.
6.96]{Kn}, says that there exists another system of simple roots which contains at most 
one non-compact root, the corresponding $r$-matrices could give rise to 
non-isomorphic Lie bialgebras, \emph{cf.} Lemma \ref{necesario} part (II) again. 
Such a possibility arises when $\invo = \jid$, or when $\invo = \jdynkin$ and
$\go$ is isomorphic to $\mathfrak{so}(2j+1, 2(n-j) - 1)$, $EI$ or $EIV$.
\end{obs}

\section{Manin triples}

In this section, we compute the Manin triples corresponding to the real absolutely simple Lie bialgebras. 
We keep the notation of the main result: $(\go, \delta)$ is an absolutely 
simple, almost factorizable, real Lie bialgebra; $\g$ is the complexification  of $\go$; 
$\invo$ is the corresponding involution, either of the form $\tdynkin$ or $\jdynkin$.

\bigbreak
We distinguish two cases:

\begin{enumerate}
\item[(a)] The bialgebra is factorizable, \emph{i.~e.} the involution $\invo$ is of the form $\tdynkin$.
\item[(b)] The bialgebra is imaginary factorizable, \emph{i.~e.} the involution $\invo$
is of the form $\jdynkin$.
\end{enumerate}

The difference between ``factorizable", case (a), and ``imaginary factorizable", case (b), 
can be read off also from the double Lie algebra: in the first case it is $\go \oplus \go$, 
in the second it is the realification $\g^{\R}$ of $\g$.

\subsection{Case (a)} In this case, the determination of the Drinfeld double and the dual 
Lie bialgebra follows from a general result from \cite{RS}. Namely, let $(\lgot, r)$ be a 
factorizable (real or complex) Lie bialgebra and let $(\quad \vert \quad)$ be the 
nondegenerate inner product on $\lgot$ induced by $r + r^{21} \in S^2 \lgot$ \cite{RS}. 
Let $r_{\pm}: \lgot^* \to \lgot$ be the maps induced by $r$, given by 
$r_{+} (\mu) = (\mu \otimes \id) r$, $r_{-} (\mu) = -(\id \otimes \mu) r$. 
The factorization map $I: \lgot^* \to \lgot$ is $I = r_+ - r_-$. 
By hypothesis, $I$ is an isomorphism and $(I (\mu)\vert I(\tau)) 
= \langle \tau, I(\mu) \rangle = \langle \mu, I(\tau) \rangle $.

\bigbreak Let $\lgot^r$ be the image of the map $\lgot^* \to \lgot \oplus \lgot$, 
$\mu \mapsto (r_{+} (\mu), r_{-} (\mu))$; it is well-known that $ r_{\pm}$ are 
Lie algebra maps, and that $\lgot^r$ is a Lie subalgebra of $\lgot \oplus \lgot$. 
Let $\diag \lgot$ be the diagonal Lie subalgebra of $\lgot \oplus \lgot$.

\begin{theorem}\label{rs-double} \cite{RS} The Manin triple corresponding to the 
factorizable Lie bialgebra $(\lgot, r)$ is $(\lgot \oplus \lgot, \diag \lgot, \lgot^r)$ 
where $\lgot \oplus \lgot$ is endowed with the bilinear form
$\langle (x,u) \vert (y,v)\rangle = (x\vert y) - (u\vert v)$, $x,y,u,v\in \lgot$. \qed
\end{theorem}

A finer description of $\lgot^r$, in terms of the Cayley transform, can be found
in \cite[Section 2.1]{Y}; see also \cite[Section 3.1]{Y} for the case of complex simple 
factorizable Lie bialgebras.

\subsection{Case (b)} We deduce from Theorem \ref{rs-double} the determination of the Manin triples 
corresponding to imaginary factorizable real Lie bialgebras.

\bigbreak
Let $(\lgo, \delta)$ be a real Lie bialgebra such that its complexification $(\lgot, \delta)$ is 
factorizable with $r\in \lgot \otimes \lgot$; let $\sigma: \lgot \to \lgot$ be the involution 
corresponding to $\lgo$. \emph{We shall assume that $\lgo$ is almost factorizable, 
see definition \ref{alm-fact}.}
Let $(\quad \vert \quad)$ be the nondegenerate inner product on $\lgot$ induced 
by $r + r^{21} = 2 i \rom \in S^2 \lgot$, and let $r_{\pm}: \lgot^* \to \lgot$ be 
the maps induced by $r$, as above. We identify $\lgo^*$ with a real subspace of $\lgot^*$; 
namely with $\{\alpha\in\lgot^*: \alpha(\lgo) \subset \R\}$. Then  
$\lgot^* = \lgo^* \oplus i\, \lgo^*$. Let $\mu \in \lgot^*$ and write 
$\mu = \alpha + i \beta$, with $\alpha, \beta \in \lgo^*$. Then 
$I(\mu) = (-\beta\otimes \id)(2\rom) + i (\alpha\otimes \id)(2\rom)$; 
in particular $\lgo = I(i\lgo^*)$. Hence, if $a = I(i\alpha), b \in \lgo$ 
then $(a\vert b) = \langle i\alpha, b\rangle \in i\R$; in other words, 
$(\quad \vert \quad)( \lgo \times \lgo) = i\R$.

\bigbreak
Consider the realification $\lgot^{\R}$ of $\lgot$. 
To avoid confusions, we denote by $x \mapsto x^{\prime}$ the multiplication 
by $i$ considered as a real linear endomorphism of $\lgor$. 
Then $\lgor = \lgo \oplus {\lgo}^{\prime}$.
%; \emph{a fortiori}, $(\lgor)^* = \lgo^* \oplus ({\lgo}^{\prime})^*$. 
The following properties are evident:
$$
x^{\prime\prime} = x, \quad [x, y^{\prime}] = [x^{\prime}, y] 
= [x,y]^{\prime}, \quad [x^{\prime}, y^{\prime}] = -[x,y], \quad \invo(x^{\prime}) 
= -\invo(x)^{\prime}, \quad x, y \in \lgor.
$$
The real bilinear form $\re(\quad \vert \quad): \lgor \times \lgor \to \R$ 
is invariant and non-degenerate; one has 
\begin{equation}\label{realpart}
2\re(u \vert v) = (u\vert v) - (\invo(u) \vert \invo(v)).
\end{equation}

\bigbreak
Let $\lgo^r := r_+(\lgo^*)$.

\begin{prop}\label{im-double} The Manin triple corresponding to the Lie bialgebra 
$(\lgo, \delta)$ is $(\lgot^{\R}, \lgo, \lgo^r)$ where $\lgot^{\R}$ is endowed 
with the bilinear form equal to $2\re(\quad \vert \quad)$. 
\end{prop}

\pf Let $\Psi : \lgot \oplus \lgot \to (\lgot^{\R})^{\Cc}$, 
$\Phi : (\lgot^{\R})^{\Cc} \to\lgot \oplus \lgot$ be given by
$$\Psi (x,y)=
\frac{1}{2}(x+\sigma (y)) + \frac{i}{2} (- x^{\prime} +\sigma (y)^{\prime}),
\qquad \Phi (u + iv) = (u + iv, \invo((u - iv)),$$
$x, y\in \lgot$, $u, v\in \lgor$. 
Notice the abuse of notation: in the argument of $\Phi$, $u+i v$ lives 
in the complexification of $\lgot^{\R}$ while in the first component 
$u+i v$ lives in $\lgot$. A straightforward computation shows that $\Phi$, 
$\Psi$ are mutually inverse isomorphisms of complex Lie algebras. We claim that
\begin{enumerate}
\item[(i)] $\Psi \diag(\lgot) = \lgo \oplus i\lgo$.
\item[(ii)] $\Psi(\lgot^r) = \lgo^r \oplus i\lgo^r$.
\item[(iii)] $\langle \Phi(u + iv) \vert \Phi(w + i z) \rangle 
= 2\re(u + iv\vert w + i z)$, for $u,v,w,z \in \lgor$, where the form 
$\langle \quad \vert \quad \rangle$ is the form defined in Theorem \ref{rs-double}, 
and $2\re(\quad \vert \quad)$ is the complexification of the real form with the same name. 
\end{enumerate}

\bigbreak
If $x\in \lgot$ then $\Psi (x, x)=\frac{1}{2}(x+\sigma (x)) 
- \frac{i}{2} (x^{\prime} +\sigma (x^{\prime}))$; this proves (i). 
Let $\mu \in \lgot^*$; say $\mu = \alpha +i \beta$, with $\alpha, \beta \in \lgo^*$. 
If we consider $r_{\pm}: \lgot^* \to \lgor$, we have 
\begin{align*}r_+(\mu) 
&= \left[(\alpha \otimes \id) \rla -  (\beta \otimes \id) \rom\right] + 
\left[(\alpha \otimes \id) \rom +  (\beta \otimes \id) \rla\right]^{\prime}, \\
r_+(\mu) ^{\prime} 
&= \left[(\alpha \otimes \id) \rla -  (\beta \otimes \id) \rom\right] ^{\prime} 
-\left[(\alpha \otimes \id) \rom +  (\beta \otimes \id) \rla\right], \\
r_-(\mu) 
& = \left[(\alpha \otimes \id) \rla +  (\beta \otimes \id) \rom\right] + 
\left[-(\alpha \otimes \id) \rom +  (\beta \otimes \id) \rla\right] ^{\prime}, \\
\invo(r_-(\mu)) 
& = \left[(\alpha \otimes \id) \rla +  (\beta \otimes \id) \rom\right] + 
\left[(\alpha \otimes \id) \rom -  (\beta \otimes \id) \rla\right] ^{\prime},\\ \invo(r_-(\mu)) ^{\prime}
& = \left[(\alpha \otimes \id) \rla +  (\beta \otimes \id) \rom\right] ^{\prime} 
+ \left[-(\alpha \otimes \id) \rom +  (\beta \otimes \id) \rla\right].
\end{align*}
Hence
$\Psi(r_+(\mu), r_+(\mu)) = 
\frac{1}{2}( r_+(\mu) +\sigma (r_-(\mu))) + \frac{i}{2} (-r_+(\mu)^{\prime} 
+\sigma (r_-(\mu))^{\prime}) = r_+(\alpha) + i r_+(\beta)$; this proves (ii).
Finally, the verification of (iii) is a straightforward computation using \eqref{realpart}. 

We conclude from the claim, by Theorem \ref{rs-double}, 
that $(\lgot^{\R}, \lgo, \lgo^r)$ 
is a Manin triple. The induced cobracket on $\lgo$ is well the initial one, again by 
Theorem \ref{rs-double}; the proof is finished.  \epf

Given a Manin triple $(\p, \p_1, \p_2)$, the isotropic Lie subalgebra $\p_2$ is not determined by
$\p_1$. Compare Proposition \ref{im-double} with the Manin triple of a compact Lie algebra
(with trivial BD-triple) in \cite{LW, M}.

\end{document}